\documentclass[12pt,a4paper]{article}
\usepackage[utf8x]{inputenc}
\usepackage[T1]{fontenc}
\usepackage{mathrsfs}
\usepackage{amsfonts}

\usepackage{soul}
\usepackage{ulem}
\usepackage{cancel}

\usepackage{comment}
\setstcolor{red}
\usepackage{xcolor}
\newcommand{\mg}{\color{magenta}}%SM
\renewcommand{\mg}{}%SM
\newcommand{\MG}{\color{magenta}}%SM
\renewcommand{\MG}{}%SM

\newcommand{\bl}{\color{blue}}%SM
\renewcommand{\bl}{}

\newcommand{\vl}{\color[RGB]{148,0,211}}%SM%violet colour
\renewcommand{\vl}{}

\newcommand{\mycomment}[1]{}
\newcommand{\R}{{\mathbb R}}

\usepackage{amsmath}
\usepackage{amsthm} %SM
\numberwithin{equation}{section}
\allowdisplaybreaks

\newtheorem{theorem}{Theorem}[section] %SM

\newtheorem{corollary}[theorem]{Corollary}
\newtheorem{lemma}[theorem]{Lemma}
\newtheorem{Rem}[theorem]{Remark} %SM
\newtheorem{defn}[theorem]{Definition} %SM
\usepackage{geometry}
\usepackage{graphicx}
\usepackage{epstopdf}
\usepackage{amsmath}
\usepackage{amssymb}
\usepackage{marginnote}
\usepackage{caption}
\usepackage{subcaption}
\usepackage{graphicx}
\usepackage{float}
\usepackage{a4wide}
\usepackage{tikz}
\usetikzlibrary{calc}
\usepackage[numbers]{natbib}
\usepackage{hyperref}
\hypersetup{
	colorlinks=true,
	linkcolor=blue,
	citecolor=blue,    
	urlcolor=blue,
} 
\usepackage{authblk}
\title{\vl Evolution Parametrix and Parametrix-Based Volume and Layer Potentials  for the Nonstationary   Variable-Coefficient Heat Equation}
\author[1]{Tamirat T. Dufera}
\author[2]{Sergey E. Mikhailov\thanks{Corresponding author: sergey.mikhailov@brunel.ac.uk}}
\author[1]{Rida T. Redi}
\affil[1]{Department of Applied Mathematics, Adama Science and Technology University, Adama, Ethiopia}
\affil[2]{Department of Mathematics, Brunel University of London, UK}
\date{November 27, 2026}
\begin{document}	
	\maketitle
	\begin{abstract}\vl
		We construct a space-time parametrix (Levi function) for a linear nonstationary heat equation with a variable coefficients. First, we prove that the associated remainder distribution can be identified with a locally integrable function. Then, we establish that the parametrix-based potentials for the variable-coefficient heat equation, namely, the initial, Newtonian, single-layer, and double-layer potentials can be expressed in terms of their well-known counterparts for the constant-coefficient heat equation. Such relations give a possibility for a subsequent analysis of the mapping properties of parametrix-based potentials. 
	\end{abstract}
	{{\bf Keywords:} \vl nonstationary heat equation; variable-coefficient; space–time parametrix;\\ 
parametrix-based volume and layer potentials.}
\\[1ex]
{\bf MSC classification:} 35K05, 35A07, 45F15
	
	\section{Introduction}
	The boundary integral equation (BIE) method is a powerful tool for solving linear partial differential equations (PDEs) with constant coefficients \citep{mclean2000, steinbach-2008, costabel1990boundary}. This method relies heavily on the availability of an explicit fundamental solution \citep{costabel1990boundary, DohrSteinbachNiino+2019+1+60}. For PDEs with variable coefficients, a fundamental solution is often unknown or cannot be {\vl constructed in an analytic form,} \citep{Evans2010, Friedman1983}. In this case, {\vl a parametrix (Levi function), see, e.g., \citep{levi1907sulle, Friedman1983},  can be used as a substitute of the fundamental solution 
 in reducing the boundary value problems to}
boundary-domain integral equations (BDIEs) \citep{ chkadua2009analysis2, mikhailov2002localized, mikhailov2006analysis, mikhailov2002localized, mclean2000}.
	
	The use of a parametrix in this way is well-known for stationary (elliptic) problems with variable coefficients. 
The appropriate analysis of the parametrix-based volume and layer potentials creates the necessary foundation for proving the equivalence of the considered boundary value problem to a system of BDIEs and for investigating its unique solvability. 

In this paper we extend such a framework, well-developed for stationary problems, see, e.g. \citep{chkadua2009analysis2, 
ayele2019analysis, mikhailov2002localized, dufera2019boundary}, to a space-time parametrix for the linear non-stationary heat equation with  a variable coefficient:
	\begin{align}\label{non-stationaryoperator}
\frac{\partial u(x, t)}{\partial t} -{\vl A{u}(x,t)}
= f(x, t), \quad \text{for } t>0,\,\, x \in \mathbb{R}^d,
	\end{align}
	where 
\begin{equation}\label{operator A}
		A{u}(x,t) :=  \sum_{i=1}^{d} \frac{\partial}{\partial x_i} \Bigl(a(x,t) \frac{\partial}{\partial x_i}u\Bigr)
	\end{equation}
and the coefficient \(a(x,t)\) 
is a given function
\citep{Evans2010, LionsMagenes1972}.
	
	The following are the main objectives of this paper.
	\begin{enumerate}
		\item To construct an explicit space-time parametrix for the operator in \eqref{non-stationaryoperator}.
		\item To prove that the corresponding remainder distribution can be identified with a locally integrable function.
		\item To establish the relation between the parametrix-based layer and volume potentials for our variable-coefficient problem with the potentials for the constant-coefficient heat equation. These include the initial, Newtonian, single-layer, and double-layer potentials.
	\end{enumerate}
	
	The results we present here provide an essential analytical foundation for the subsequent study of equivalence of the considered variable-coefficient initial-boundary value problem to some systems of space-time boundary-domain integral equations, and for the investigation of their unique solvability, in another paper.

		\section{A space-time parametrix and remainder 
}\label{S2}

		We will construct a space-time parametrix for the linear non-stationary diffusion problem with variable coefficients \eqref{non-stationaryoperator} by the method of ``frozen coefficient.''
		
We assume that $a(x,t)$ 	in \eqref{non-stationaryoperator} is  a positive and  bounded  scalar function, $a\in C^\infty({\R^{d+1}})$.
and furthermore, 
$$\vl 0<a_{\min}\le a(x,t)\le a_{\max} < \infty.$$  
Then, {\vl denoting $a_{ij}:=\delta_{ij}a$},
we have the uniform parabolicity of equation \eqref{non-stationaryoperator}:	
		\begin{equation}\label{m and M}	
			a_{\min}|\eta|^2 \leq \sum_{i, j=1}^d a_{i j}(x,t) \eta_i \eta_j 
{\vl=a(x,t)\sum_{i}^d \eta_i^2}
\leq a_{\max}|\eta|^2, \quad  x, \eta \in \mathbb{R}^{d}.	
		\end{equation}

Let us introduce the differential operator 
				\begin{equation}\label{3.1}	
{\vl \mathcal{L}}=	\mathcal{L}_{x,t} := \dfrac{\partial}{\partial t}-A
=\dfrac{\partial}{\partial t} - \sum_{i=1}^{d} \frac{\partial}{\partial x_i} \Bigl(a(x,t) \frac{\partial}{\partial x_i}\Bigr).
		\end{equation}
Let us also define the operator formally adjoint to $\mathcal{L}$, 
	\begin{equation}\label{adjoint of L}
		\mathcal{L}^* =	\mathcal{L}^*_{x,t} := -\partial_t - A
=-\dfrac{\partial}{\partial t} - \sum_{i=1}^{d} \frac{\partial}{\partial x_i} \Bigl(a(x,t) \frac{\partial}{\partial x_i}\Bigr).
	\end{equation}
		Analogous to the definition of a parametrix for linear elliptic PDEs, we define  space-time parametrices of the operators $\mathcal{L}$ and $\mathcal{L}^*$ as follows.

		\begin{defn}\label{definition of parametrix}\ 
(i) A  parametrix of the differential operator $\mathcal{L}$ in \eqref{3.1}, is for any fixed $(y,\tau)\in \R^d\times\R$ a distribution $P(x,y,t,\tau)$  in  $\mathbb{R}^{d+1}$ satisfying
			\begin{equation}\label{para1}
				\mathcal{L}_{x,t} P(x,y,t,\tau) = \delta(x - y, t - \tau) + R(x,y,t,\tau), 
			\end{equation}  
			where $ \delta(\cdot, \cdot)$ is the Dirac-delta distribution in $\mathbb{R}^{d+1}$. 
{\vl Moreover, for any fixed $(y,\tau)\in \R^d\times\R$ the remainder $R(x,y,t,\tau)$  is a locally  integrable function of $(x,t)\in \R^d\times\R$}.

(ii)	A parametrix of a differential operator $\mathcal{L}^*$ in \eqref{adjoint of L} is for any fixed $(y,\tau)\in \R^d\times\R$ a distribution $P^*(x,y,t,\tau)$  in  $\mathbb{R}^{d+1}$ satisfying
			\begin{equation}\label{para1*}
				\mathcal{L}^*_{x,t} P^*(x,y,t,\tau) = \delta(x - y, t - \tau) + R^*(x,y,t,\tau).
			\end{equation}  
{\vl Moreover, for any fixed $(y,\tau)\in \R^d\times\R$ the remainder $R^*(x,y,t,\tau)$  is a locally  integrable function of $(x,t)\in \R^d\times\R$}.
		\end{defn}

To construct a parametrix for the operator  \eqref{non-stationaryoperator}, we use the well known  fundamental solution for the partial differential operator $\mathcal L_{\underline a}$, with a constant coefficient $\underline{a}$ 
		\begin{equation}\label{4}
			\mathcal L_{\underline a}= \frac{\partial }{\partial t}- \underline{a}\sum_{i=1}^{d}  \frac{\partial^2 }{\partial x^2_{i}},
		\end{equation}                       
(see, e.g., \cite[Section 10.6, Eq. (25)]{Vladimirov1971}, \cite[Chapter 1, Eq. (0.2)]{Friedman1983}) given by
		\begin{equation} \label{5.1}
			P_{\underline{a}}(x, t)=\bar P(x, t; \underline{a}) = \frac{\Theta(t)}{(4\pi \underline{a} t)^{d/2}}\exp \left(\frac{-|x|^{2}}{4\underline{a}t}\right),\,\, x \in \mathbb{R}^d,\ t\in \mathbb R 
		\end{equation}
		for $d\ge 1$, that is,
\begin{equation}\label{fundamental solution acts to La 0}
			\mathcal L_{\underline{a}}P_{\underline{a}}(x,t)=\delta(x,t).
		\end{equation}
Here $\Theta$ is the Heaviside step function,
$
			\Theta(t) = \begin{cases} 1 , & t > 0 , \\ 0,& \le  0  \end{cases}
$ .

The function 
\begin{equation} \label{9.1}
P_{\underline{a}}(x,y, t,\tau)
=\bar P(x-y, t-\tau;\underline{a}) 
= \frac{\Theta(t-\tau)}{(4\pi \underline{a}\ (t - \tau))^{d/2}}\exp \left(\frac{-|x - y|^{2}}{4\underline{a}\ (t - \tau)}\right),
\, \,x,y \in \mathbb{R}^d, \ \ t ,\tau\in\mathbb R  
		\end{equation}
is evidently a (shifted) fundamental solution for the partial differential operator $\mathcal L_{\underline a}$, with the right-hand side delta distribution having support at $(y,\tau)$, i.e.,
		\begin{equation}\label{fundamental solution acts to La}
			\underline{\mathcal L}_{\underline a;x,t}P_{\underline{a}}(x,y, t,\tau)=\delta(x-y,t-\tau).
		\end{equation}

Note also that the function 
\begin{align} \label{9.1*}
			P^*_{\underline{a}}(x,y, t,\tau) 
&=\bar P(x-y,\tau-t;\underline{a}) 
\nonumber\\
&= \frac{\Theta(\tau-t)}{(4\pi \underline{a}\ (\tau-t))^{d/2}}\exp \left(\frac{-|x - y|^{2}}{4\underline{a}\ (\tau-t)}\right),
\, \,x,y \in \mathbb{R}^d, \ \ t ,\tau\in\mathbb R  
		\end{align}
is a (shifted) fundamental solution for the partial differential operator $\mathcal L^*_{\underline a}$, with the right-hand side delta distribution having support at $(y,\tau)$, i.e.,
		\begin{equation}\label{fundamental solution acts to La*}
			\underline{\mathcal L}^*_{\underline a;x,t}P^*_{\underline{a}}(x,y, t,\tau)=\delta(x-y,t-\tau)
		\end{equation}
where, as in \eqref{adjoint of L}, 
\begin{equation*}\label{4*}
			\mathcal L^*_{\underline a}= -\frac{\partial }{\partial t}- \underline{a}\sum_{i=1}^{d}  \frac{\partial^2 }{\partial x^2_{i}}.
		\end{equation*}  

		Now we construct a parametrix of the variable-coefficient differential operator  $\mathcal{L}$ in \eqref{3.1}. 
This can be achieved by the method of ``frozen coefficient''  (cf., e.g.,\cite{levi1907sulle}, \cite{Friedman1983}, \cite{mikhailov2002localized}, \cite{chkadua2009analysis2}). 
That is, replacing the constant $\underline{a}$ in \eqref{9.1} by $a(y,\tau)$, we define 
	\begin{align} \label{12.2}
				P(x,y,t,\tau)
&=\bar P(x-y,t-\tau;a(y,\tau)) 
\nonumber\\
&= \frac{\Theta(t-\tau)}{(4\pi a(y,\tau) (t - \tau))^{d/2}}\exp \left(\frac{-|x - y|^{2}}{4a(y,\tau) (t - \tau)}\right), 
\, \,x,y \in \mathbb{R}^d, \ \ t ,\tau\in\mathbb R.
			\end{align}

Let us set 
\begin{align}\label{chvar}
z = x - y,\quad \zeta = t -\tau.
\end{align} 
For $\zeta> 0$, let us also
perform the change of variables 
\begin{align}\label{E2.12}
w = z/\sqrt{\zeta}, \text{ hence } dx =dz= \zeta^{d/2}dw.
\end{align}
Then relation \eqref{12.2} leads to the following estimate for $\zeta>0$,
			\begin{equation}\label{Pest}
				|P(x,y,t,\tau)|
= (4\pi a(y,\tau)\zeta)^{-d/2} e^{-|w|^2/(4a(y,\tau))}
\le (4\pi a_{\min}\zeta)^{-d/2} e^{-|w|^2/(4a_{\max})} .
			\end{equation}

Let us recall the well-known result obtained for any $m>-d$ and $\alpha>0$ by using the spherical coordinates. 	
\begin{align}\label{Gamma}
\int_{\mathbb{R}^d}|w|^{m}e^{-\alpha|w|^2}\,dw
= |S_{d-1}|\int_0^\infty r^{m+d-1}e^{-\alpha r^2}\,dr	
=\frac{1}{2}|S_{d-1}| \alpha^{-(m+d)/2}\Gamma\left(\frac{m+d}{2}\right),
\end{align}
where  $\Gamma$ if the Euler Gamma-function and 
$
|S_{d-1}|=\frac{2\pi^{d/2}}{\Gamma(d/2)}
$ 
is surface area of the unit sphere, $S_{d-1}$.

\begin{Rem} The function $P(x,y,t,\tau)$ has a singularity at $x=y$ {\vl as $t\searrow \tau$}, which is weak since for any fixed $(y,\tau)\in\R^{d+1}$, the function $P(x,y,t,\tau)$ defined by \eqref{12.2} is locally integrable in $(x,t)\in\R^{d+1}$.
Indeed, 
for some $T_0>0$, let us take the integration domain  $\R^d\times(0,T_0)$ in $(z,\zeta)$. 
Then relation \eqref{12.2} leads to the estimate
\begin{multline*}
\int_{\tau}^{\tau+T_0}\int_{\R^d}|P(x,y,t,\tau)|dxdt
=	\int_{0}^{T_0}\int_{\R^d}\bar P(z,\zeta;a(y,\tau))dzd\zeta 
\\
=	\int_{0}^{T_0}\int_{\R^d} (4\pi a(y,\tau)\zeta)^{-d/2} e^{-|w|^2/(4a(y,\tau))} \zeta^{d/2}dwd\zeta
\\
\le (4\pi a_{\min})^{-d/2} \int_{0}^{T_0}\int_{\R^d}e^{-|w|^2/(4a_{\max})} dwd\zeta
\\
=\frac12\left(\frac{a_{max}}{\pi a_{\min}}\right)^{d/2} {T_0} |S_{d-1}|\Gamma({d}/{2})
=\left(\frac{a_{max}}{a_{\min}}\right)^{d/2} {T_0} ,
			\end{multline*}
where we employed \eqref{Gamma} with $m=0$.
This means $P(x,y,t,\tau)$ is locally integrable.
\end{Rem}

		\begin{theorem}\label{T2.2}
		\vl For any fixed $(y,\tau)\in\R^{d+1}$, the function $P(x,y,t,\tau)$ defined by \eqref{12.2} {\vl is a parametrix for the differential operator $\mathcal{L}_{x,t}$. 
It satisfies equation \eqref{para1} in the sense of distributions, with the remainder}
\begin{align}\label{Reminder0}
				R(x,y,t,\tau)=\sum_{i=1}^{d}\dfrac{\partial  }{\partial x_i}  \left([a(y,\tau) -a(x,t)]\dfrac{\partial  P(x,y,t,\tau)}{\partial x_i}\right)
\end{align}
and the distribution $R$ can be identified with the locally integrable function of $(x,t)\in\R^{d+1}$.
Moreover,
\begin{align}\label{Reminder}
	R(x,y,t,\tau)&=\sum_{i=1}^{d}  -\dfrac{\partial a(x,t)}{\partial x_i} \dfrac{\partial  P(x,y,t,\tau)}{\partial x_i}+\widetilde R(x,y,t,\tau),
\\
\label{Reminder~}
\text{where}\quad\widetilde R(x,y,t,\tau)&:=[a(y,\tau) -a(x,t)]\Delta_x P(x,y,t,\tau).
\end{align}
Explicitly,
\begin{multline}\label{Reminder-e}
R(x,y,t,\tau)= \frac{\Theta(t-\tau)}{(4\pi a(y,\tau) (t - \tau))^{d/2}}\exp \left(\frac{-|x - y|^{2}}{4a(y,\tau) (t - \tau)}\right)
\Bigg(\sum_{i=1}^{d}  \dfrac{\partial  a(x,t)}{\partial x_i} \dfrac{x_i-y_i  }{2a(y,\tau)(t-\tau)}
\\
+[a(y,\tau) -a(x,t)] \left(-\frac{d}{2a(y,\tau)(t-\tau)} + \frac{|x-y|^2}{4a(y,\tau)^2(t-\tau)^2}\right)\Bigg).
			\end{multline}
		\end{theorem}
		\begin{proof}
Let us first show that the function $P=P(x,y,t,\tau)$ 
satisfies equation 	\eqref{para1} in the sense of distributions.		
			Indeed,
			\begin{align}\label{14.11}
				&\mathcal{L}_{x,t}P(x,y,t,\tau)=
				\frac{\partial P(x,y,t,\tau)}{\partial t} -\sum_{i=1}^{d}\frac{\partial}{\partial x_{i}}\left(a(x,t) \frac{\partial P(x,y,t,\tau)}{\partial x_{i}}\right)  
\nonumber\\
					&= \left( \dfrac{\partial }{\partial t} - a(y,\tau)\Delta_x\right)P(x,y,t,\tau) + \sum_{i=1}^{d}\dfrac{\partial}{\partial x_i}\left([a(y,\tau) -a(x,t)]\dfrac{\partial P(x,y,t,\tau)}{\partial x_i}\right)
\nonumber\\ 
					&= \delta(x-y, t-\tau) + R(x,y,t,\tau),
			\end{align}
			where 	
	\begin{align}\label{E2.13}
	R(x,y,t,\tau)&=\sum_{i=1}^{d}\dfrac{\partial  }{\partial x_i}  \left([a(y,\tau) -a(x,t)]\dfrac{\partial  P(x,y,t,\tau)}{\partial x_i}\right).
\end{align}
	In the last line of \eqref{14.11}, we have used the relation $\left(\partial_t-a(y,\tau)\Delta_x \right)P=\delta(x-y,t-\tau)$, see \eqref{fundamental solution acts to La}. 

For $P(x,y,t,\tau)$  to be a parametrix,
we have to show that the remainder $R(x,y,t,\tau)$ can be identified with a locally integrable (weakly singular) function in the space-time variables $(x,t)\in\R^{d+1}$. 

Let $\phi \in \mathcal{D}(\R^{d+1})\vl:=C^\infty_c(\R^{d+1})$ be an arbitrary test function. 
Let us fix $(y,\tau)\in\R^{d+1}$ and consider $R(x,y,t,\tau)$ as a distribution in $(x,t)\in\R^{d+1}$.
Then, slightly abusing the notations, we have
\begin{multline}\label{Rphi}
\langle R(\cdot, y, \cdot, \tau),\phi\rangle_{\R^{d+1}}
=\left\langle \sum_{i=1}^{d}\dfrac{\partial  }{\partial x_i}  \left([a(y,\tau) -a(x,t)]\dfrac{\partial  P(x,y,t,\tau)}{\partial x_i}\right),
\phi(x,t)\right\rangle_{(x,t)\in\R^{d+1}}
\\
\shoveright{=\left\langle  P(x,y,t,\tau),
\sum_{i=1}^{d}\dfrac{\partial  }{\partial x_i}\left([a(y,\tau) -a(x,t)]\dfrac{\partial  }{\partial x_i} \phi(x,t)\right)\right\rangle_{(x,t)\in\R^{d+1}} }
\\
=\int_{\R^{d+1}}  P(x,y,t,\tau)
\sum_{i=1}^{d}\dfrac{\partial  }{\partial x_i}\left([a(y,\tau) -a(x,t)]\dfrac{\partial  }{\partial x_i} \phi(x,t)\right)dxdt.\ 
\end{multline}
The last equality follows from  the inclusion $\phi \in \mathcal{D}(\R^{d+1})$ and the local integrability of $P$.

Let $B_\epsilon$ be the ball of radius $\epsilon$ for the variable $w$ in $\R^d$ centred in the origin, and $n(w)$ be the inward unit normal to the ball boundary $\partial B_\epsilon$. Then $B_{\epsilon\sqrt{\zeta}}$ is the corresponding ball of radius $\epsilon\sqrt{\zeta}$ for the variable $z$ and $n(z)=n(w)$ is the inward unit normal to the ball boundary $\partial B_{\epsilon\sqrt{\zeta}}$.
Transformations \eqref{Rphi} can be continued as follows,
\begin{multline}\label{Rphi2}
\langle R(\cdot, y, \cdot, \tau),\phi\rangle_{\R^{d+1}}
\\
=\lim_{\epsilon\to 0}\int_{\R}\int_{\R^{d}\setminus B_{\epsilon\sqrt{\zeta}}} \sum_{i=1}^{d} \bar P(z,\zeta;a(y,\tau))
\dfrac{\partial  }{\partial z_i}\left([a(y,\tau)-a(z+y,\zeta+\tau)]\dfrac{\partial  }{\partial z_i} \phi(z+y,\zeta+\tau)\right)dzd\zeta
\\
\shoveleft
{=\lim_{\epsilon\to 0}\int_{0}^\infty\int_{\partial B_{\epsilon\sqrt{\zeta}}} \sum_{i=1}^{d} \bar P(z,\zeta;a(y,\tau))
n_i(z)[a(y,\tau)-a(z+y,\zeta+\tau)]\dfrac{\partial  }{\partial z_i} \phi(z+y,\zeta+\tau)d\sigma_zd\zeta}
\\-\lim_{\epsilon\to 0}\int_{0}^\infty\int_{\R^{d}\setminus B_{\epsilon\sqrt{\zeta}}}  \sum_{i=1}^{d}\left(\dfrac{\partial  }{\partial z_i}\bar P(z,\zeta;a(y,\tau))\right)
[a(y,\tau)-a(z+y,\zeta+\tau)]\dfrac{\partial  }{\partial z_i} \phi(z+y,\zeta+\tau)dzd\zeta
\\
\shoveleft{=\lim_{\epsilon\to 0}\int_{0}^\infty\int_{\partial B_{\epsilon\sqrt{\zeta}}} \sum_{i=1}^{d} \bar P(z,\zeta;a(y,\tau))
n_i(z)[a(y,\tau)-a(z+y,\zeta+\tau)]\dfrac{\partial  }{\partial z_i} \phi(z+y,\zeta+\tau)d\sigma_zd\zeta}
\\
{-\lim_{\epsilon\to 0}\int_{0}^\infty\int_{\partial B_{\epsilon\sqrt{\zeta}}} \sum_{i=1}^{d}\left(\dfrac{\partial  }{\partial z_i}\bar P(z,\zeta;a(y,\tau))\right)
n_i(z)[a(y,\tau)-a(z+y,\zeta+\tau)] \phi(z+y,\zeta+\tau)d\sigma_zd\zeta}
\\
\shoveright{+\lim_{\epsilon\to 0}\int_{0}^\infty\int_{\R^{d}\setminus B_{\epsilon\sqrt{\zeta}}}  \sum_{i=1}^{d}\dfrac{\partial  }{\partial z_i}\left([a(y,\tau)-a(z+y,\zeta+\tau)]
\dfrac{\partial  }{\partial z_i}\bar P(z,\zeta;a(y,\tau))\right)
 \phi(z+y,\zeta+\tau)d\sigma_zd\zeta}
\\
=I_1+I_2+I_3,
\end{multline}
where
\begin{multline}\label{I1}
I_1:=\lim_{\epsilon\to 0}\int_{0}^\infty\int_{\partial B_{\epsilon}}\Big[ \sum_{i=1}^{d} \bar P(z,\zeta;a(y,\tau))
n_i(z)[a(y,\tau)-a(z+y,\zeta+\tau)]
\\
\times \dfrac{\partial  }{\partial z_i}\phi(z+y,\zeta+\tau)\Big]_{z=w\sqrt{\zeta}}\zeta^{(d-1)/2}d\sigma_wd\zeta,
\end{multline}
\begin{multline}\label{I2}
I_2:=-\lim_{\epsilon\to 0}\int_{0}^\infty\int_{\partial B_{\epsilon}}\Big[ \sum_{i=1}^{d}\left(\dfrac{\partial  }{\partial z_i}\bar P(z,\zeta;a(y,\tau))\right)
n_i(z)[a(y,\tau)-a(z+y,\zeta+\tau)]
\\
\times \phi(z+y,\zeta+\tau)\Big]_{z=w\sqrt{\zeta}}\zeta^{(d-1)/2}d\sigma_w d\zeta,
\end{multline}
\begin{multline}\label{I3}
I_3:
=\lim_{\epsilon\to 0}\int_{0}^\infty\int_{\R^{d}\setminus B_{\epsilon}}
 \Big[ \sum_{i=1}^{d}\dfrac{\partial  }{\partial z_i}\left([a(y,\tau)-a(z+y,\zeta+\tau)]
\dfrac{\partial  }{\partial z_i}\bar P(z,\zeta;a(y,\tau))\right)
\\
\shoveright{\times\phi(z+y,\zeta+\tau)\Big]_{z=w\sqrt{\zeta}}\zeta^{d/2}dwd\zeta}
\\
=\lim_{\epsilon\to 0}\int_{0}^\infty\int_{\R^{d}\setminus B_{\epsilon}} \Big[R(z+y,y,\zeta+\tau,\tau)
 \phi(z+y,\zeta+\tau)\Big]_{z=w\sqrt{\zeta}}\zeta^{d/2}dwd\zeta.
\end{multline}
The last equality in the definition of $I_3$ follows from \eqref{E2.13}.

Let ${\rm supp}(\phi)\subset B_{R_\phi}\times(-T_\phi,T_\phi)$ for some $R_\phi,T_\phi>0$.
Let  $|\phi(x,t)| \le M_\phi$, $|\nabla_x \phi(x,t|\le  M'_\phi$, $|\nabla_x a(x,t)|,|\partial_t  a(x,t)|\le M'_a$, and  hence
			\begin{equation}\label{E2.14}
				|a(y,\tau) - a(x,t)| \le M'_a(|z|+\zeta)=M'_a|w\sqrt{\zeta}|+M'_a\zeta
			\end{equation}
for $(x,t)\in B_{R_\phi}\times(-T_\phi,T_\phi)$ and $\zeta>0$.

Denoting $T_{\phi,\tau}:=\max\{T_\phi-\tau,0\}$, due to \eqref{I1} and estimate \eqref{Pest}, we have
\begin{multline}\label{E2.15}
|I_1|\le
\lim_{\epsilon\to 0}\int_{0}^\infty\int_{\partial B_{\epsilon}} (4\pi a_{\min}\zeta)^{-d/2} e^{-\epsilon^2/(4a_{\max})} M'_a
(\epsilon\sqrt{\zeta}+\zeta)
\\
\shoveright{\times\Big|\nabla_z\phi(z+y,\zeta+\tau)\Big|_{z=w\sqrt{\zeta}}\zeta^{(d-1)/2}d\sigma_wd\zeta}
\\
\le\lim_{\epsilon\to 0} \int_0^{T_{\phi,\tau}}\int_{\partial B_{\epsilon}} (4\pi a_{\min})^{-d/2} M'_a(\epsilon\sqrt{\zeta}+\zeta)
 M'_\phi \zeta^{-1/2}d\sigma_w d\zeta
\\
=\lim_{\epsilon\to 0}|S_{d-1}|\epsilon^{d-1} (4\pi a_{\min})^{-d/2}  M'_a(\epsilon T_{\phi,\tau}+\frac23 T_{\phi,\tau}^{3/2}) M'_\phi =0
\end{multline}

For spatial derivatives of the parametrix at $t>\tau$, we have from \eqref{12.2}, 
\begin{align}\label{nablaP}
				\nabla_x P(x,y,t,\tau) &= -\frac{x-y}{2a(y,\tau)(t-\tau)}P(x,y,t,\tau), 
				\\
\label{deltaP}
				\Delta_x P(x,y,t,\tau) &= \left(-\frac{d}{2a(y,\tau)(t-\tau)} + \frac{|x-y|^2}{4a(y,\tau)^2(t-\tau)^2}\right)P(x,y,t,\tau).
\end{align}

By \eqref{nablaP} and \eqref{Pest}, for $\zeta>0$ we obtain
			\begin{equation}
				|\nabla_z P|
\le  \left|\frac{w}{2a(y,\tau)\sqrt{\zeta}}\right| (4\pi a_{\min}\zeta)^{-d/2} e^{-|w|^2/(4a_{\max})} .
			\end{equation}
Then employing \eqref{E2.14} to \eqref{I2}, we obtain
\begin{multline}\label{E2.16}
|I_2|
\le\lim_{\epsilon\to 0}\int_{0}^\infty\int_{\partial B_{\epsilon}}\frac{\epsilon}{2a_{\min}\sqrt{\zeta}} (4\pi a_{\min}\zeta)^{-d/2} e^{-\epsilon^2/(4a_{\max})} M'_a(\epsilon\sqrt{\zeta}+\zeta)
\\
\shoveright{\times \Big|\phi(z+y,\zeta+\tau)\Big|_{z=w\sqrt{\zeta}}\zeta^{(d-1)/2}d\sigma_wd\zeta}
\\
\le\lim_{\epsilon\to 0} \int_0^{T_{\phi,\tau}}\int_{\partial B_{\epsilon}}\frac{\epsilon}{2a_{\min}} (4\pi a_{\min})^{-d/2} M'_a
(\epsilon\zeta^{-1/2}+1) M_\phi d\sigma_w d\zeta
\\
=\lim_{\epsilon\to 0}|S_{d-1}|\epsilon^{d-1}\frac{\epsilon}{2a_{\min}} (4\pi a_{\min})^{-d/2}  M'_a 
(\epsilon 2\sqrt{T_{\phi,\tau}}+T_{\phi,\tau})M_\phi =0.
\end{multline}

Further, from  \eqref{deltaP} we get 
\begin{equation}
				|\Delta_z P| \le \left(\frac{d}{2a(y,\tau)\zeta} + \frac{|z|^2}{4a(y,\tau)^2\zeta^2}\right)P.
			\end{equation}
and hence, due to \eqref{E2.14} and \eqref{nablaP},
			\begin{align}
				|R(z+y,y,\zeta+\tau,\tau)|
				&\le |a(y,\tau)-a(z+y,\zeta+\tau)|\,|\Delta_z P| + |\nabla_z a(z+y,\zeta+\tau)|\,|\nabla_z P| \nonumber\\
				&\le M'_a(|z|+\zeta)\left(\frac{d}{2a_{\min}\zeta} + \frac{|z|^2}{4a_{\min}^2\zeta^2}\right)P 
				+ M'_a\frac{|z|}{2a_{\min}\zeta}P \nonumber\\
&= M'_a\!\left[\frac{d+1}{2a_{\min}}\frac{|z|}{\zeta} 	+ \frac{1}{4a_{\min}^2}\frac{|z|^{3}}{\zeta^2}
+\frac{d}{2a_{\min}} + \frac{|z|^2}{4a_{\min}^2\zeta}\right]P.
			\end{align}
for $\zeta>0$.
By \eqref{E2.12} and \eqref{Pest}, this implies  
\begin{multline}\label{E2.21}
\int_{0}^\infty\int_{\R^{d}\setminus B_{\epsilon}} 
\Big|R(z+y,y,\zeta+\tau,\tau) \phi(z+y,\zeta+\tau)\Big|_{z=w\sqrt{\zeta}}\zeta^{d/2}dwd\zeta
\\
\shoveleft{\le\int_{0}^{T_{\phi,\tau}}\int_{\R^{d}\setminus B_{\epsilon}} 
M'_a\!\left[\frac{d+1}{2a_{\min}}\frac{|w|}{\sqrt{\zeta}} 	+ \frac{1}{4a_{\min}^2}\frac{|w|^{3}}{\sqrt{\zeta}}
+\frac{d}{2a_{\min}} + \frac{|w|^2}{4a_{\min}^2}\right]}
\\
\shoveright{\times 
(4\pi a_{\min}\zeta)^{-d/2} e^{-|w|^2/(4a_{\max})}\Big|\phi(z+y,\zeta+\tau)\Big|_{z=w\sqrt{\zeta}}\zeta^{d/2}dwd\zeta}
\\
\shoveleft{
\le\int_{0}^{T_{\phi,\tau}}\int_{\R^{d}\setminus B_{\epsilon}} 
M'_a\!\left[\frac{d+1}{2a_{\min}}\frac{|w|}{\sqrt{\zeta}} 	+ \frac{1}{4a_{\min}^2}\frac{|w|^{3}}{\sqrt{\zeta}}
+\frac{d}{2a_{\min}} + \frac{|w|^2}{4a_{\min}^2}\right]
}
\\
\shoveright{\times 
(4\pi a_{\min})^{-d/2} e^{-|w|^2/(4a_{\max})}M_\phi dwd\zeta
}
\\
\shoveleft{
=\int_{\R^{d}\setminus B_{\epsilon}} 
M'_a\!\left[\frac{d+1}{a_{\min}}{|w|}{\sqrt{T_{\phi,\tau}}} 
+ \frac{1}{2a_{\min}^2}{|w|^{3}}{\sqrt{T_{\phi,\tau}}}
+\frac{d}{2a_{\min}}T_{\phi,\tau} 
+ \frac{|w|^2}{4a_{\min}^2}T_{\phi,\tau}\right]
}
\\
\shoveright{\times 
(4\pi a_{\min})^{-d/2} e^{-|w|^2/(4a_{\max})}M_\phi dw
}
\\
\le\int_{\R^{d}} 
\!\left[C_0+ C_1|w| +C_2|w|^2+ C_3|w|^3\right]
 e^{-\alpha|w|^2}
dw,
\end{multline}
where $\alpha:={1}/{(4a_{\max})}$.
Then relations \eqref{Gamma} imply that the right hand side in \eqref{E2.21} is bounded and hence the value of $I_3$ is well defined and bounded as well.
Together with relations \eqref{Rphi2} and estimates \eqref{E2.15} and \eqref{E2.16} this means that 
\begin{align}
\langle R(\cdot , y, \cdot , \tau),\phi\rangle_{\R^{d+1}}
&=I_3
=\int_{\tau}^\infty\int_{\R^{d}}  R(x,y,t,\tau)\phi(x,t)dtdx
\end{align}
and integrals in the right-hand side are well defined in the Lebesgue sense for any $\phi \in \mathcal{D}(\R^{d+1})$.

Moreover, \eqref{E2.13}, \eqref{nablaP} and \eqref{deltaP} imply \eqref{Reminder}--\eqref{Reminder-e}.
\end{proof}

			For a constant coefficient $a $, we have $ a(x,t) = a(y,\tau)=\underline{a}$, which implies the differential operator $\mathcal{L}_{x,t}$ reduces to the constant-coefficient heat operator $(\partial_t - \underline{a}\Delta_x)$. 
In this case, $P(x,y,t,\tau)$ becomes its fundamental solution, and \mbox{$R(x,y,t,\tau) \equiv 0$.}

{\vl Making corresponding changes in the reasoning of Theorem \ref{T2.2}, we arrive at the following corollary.}
	\begin{corollary}
\vl
For any fixed $(y,\tau)\in\R^{d+1}$, the function 
\begin{multline}\label{E2.39}
P^*(x,y,t,\tau):
=\bar P(x-y,\tau-t;a(y,\tau)) 
\\
=\frac{\Theta(\tau-t)}{(4\pi a(y,\tau) (\tau-t))^{d/2}}\exp \left(\frac{-|x - y|^{2}}{4a(y,\tau) (\tau-t)}\right), 
\, \,x,y \in \mathbb{R}^d, \ \ t ,\tau\in\mathbb R, 
\end{multline}
is a parametrix for the differential operator 	 $\mathcal{L}^*_{x,t}$ in \eqref{adjoint of L}, and the corresponding remainder is expressed as
			\begin{align}\label{Reminder*0}
	R^*(x,y,t,\tau)
=\sum_{i=1}^{d}\dfrac{\partial  }{\partial x_i}  \left([a(y,\tau) -a(x,t)]\dfrac{\partial  P^*(x,y,t,\tau)}{\partial x_i}\right), 
\end{align}
where the distribution $R^*$ can be identified with the locally integrable function of $(x,t)\in\R^{d+1}$.
Moreover,
\begin{align}\label{Reminder*}
	R^*(x,y,t,\tau)&=\sum_{i=1}^{d} -\dfrac{\partial a(x,t)}{\partial x_i}\dfrac{\partial  P^*(x,y,t,\tau)}{\partial x_i}+\widetilde R^*(x,y,t,\tau),
\\
\label{Reminder*~}
\text{where}\quad\widetilde R^*(x,y,t,\tau)&:=[a(y,\tau) -a(x,t)]\Delta_x P^*(x,y,t,\tau).
\end{align}
Explicitly,
	\begin{multline}\label{Reminder*-e}
R^*(x,y,t,\tau)
=\frac{\Theta(\tau-t)}{(4\pi a(y,\tau) (\tau-t))^{d/2}}\exp \left(\frac{-|x - y|^{2}}{4a(y,\tau) (\tau-t)}\right)
\\
\shoveleft{\quad\times\Bigg(\sum_{i=1}^{d}  \dfrac{\partial  a(x,t)}{\partial x_i} \dfrac{x_i-y_i  }{2a(y,\tau)(\tau-t)}}
\\
+[a(y,\tau) -a(x,t)] \left(-\frac{d}{2a(y,\tau)(\tau-t)} + \frac{|x-y|^2}{4a(y,\tau)^2(\tau-t)^2}\right)\Bigg).
	\end{multline}
\end{corollary}
	
		\section{Classical Green's Identities and Parametrix-Based\\
 Potentials}

\subsection{Green's Identities}
	Let $\Omega \subset \mathbb{R}^d$ be a Lipschitz spatial bounded domain, $d\ge2$ be an integer dimension. Denote by $\partial\Omega$ the boundary of $\Omega$, and let $\overline{\Omega}$ denote the domain closure. Fix a final time $T>0$ and define the space-time cylinder $Q := \Omega \times (0,T)$ with lateral boundary $\Sigma := \partial\Omega \times (0,T)$.

For $u \in C(\overline{Q})$ (respectively $C(\overline{Q^c})$), the interior/exterior Dirichlet traces on $\partial\Omega$ are 
$$\gamma^{cl \pm} u(x,t) := \lim_{\Omega^\pm \ni \tilde x \to x \in \partial\Omega} u(\tilde x,t).$$ 
	
 {\mg For} $u \in C^1(\overline{Q})$ (respectively $C^1(\overline{Q^c})$), at points $x$ of the Lipschitz boundary $\partial\Omega$,
were the unit outward normal  $n^+(x)$ to $\partial\Omega$ does exist, the {\mg classical} interior/exterior co-normal derivatives on  $\partial\Omega$ are
	\begin{align}\label{Tcl}
	\mathcal{T}^{cl\pm} u(x,t) := \lim_{\Omega^{\pm} \ni \tilde x \to x} \sum_{i=1}^d n_i^\pm(x)\, a(\tilde x)\, \partial_{\tilde x_i} u(\tilde x,t),
	\end{align}
	so that $\mathcal{T}^{cl\pm} u = a(x)\,\partial_{n^\pm} u|_\Sigma$, where  $n^-=-n^+$.

By applying the Gauss divergence theorem,
we obtain the following classical Green’s identities.

For any $u \in C^{2}(\overline{Q})$ and \( v \in C^{1}(Q) \cap C(\overline{Q}) \), the {\MG first classical Green identity} holds
\begin{multline}\label{classicalFirstGreenIdentity}
		\int_{0}^{T}\!\!\int_{\Omega} v \,\partial_t u\, dx\,dt
		+ \int_{0}^{T}\!\!\int_{\Omega} a(x,t)\,\nabla u \cdot \nabla v \,dx\,dt
\\
		= \int_{0}^{T}\!\!\int_{\Omega} v\mathcal{L}u\,dx\,dt 
		\quad + \int_{0}^{T}\!\!\int_{\partial\Omega} v\,a(x,t)\partial_{n_x}u\,d\sigma_x\,dt ,
\end{multline}
After integration by parts in $t$, \eqref{classicalFirstGreenIdentity} gives
\begin{multline}\label{integrating-by-part}
	\int_{0}^{T}\!\!\int_{\Omega}\left(a(x)\nabla u \cdot \nabla v-u \partial_t v\right) \,d x d t
+\int_{\Omega} u(x,T) v(x,T) d x
\\
-\int_{\Omega} u(x,0) v(x,0) d x-\int_{\Sigma} v \,a(x,t)\partial_{n_x} u\, d t d \sigma
=\int_{0}^{T}\!\!\int_{\Omega} v\mathcal{L}u \,d x d t.
\end{multline}

Similarly, for any $u,v \in C^{2}(\overline{Q})$
we obtain the following classical Green’s second identity,
\begin{equation}\label{E4.13}
	\begin{aligned}
		&\int_{\Omega} u(x, T) v(x, T)\,dx - \int_{\Omega} u(x, 0) v(x, 0)\,dx \\
		&\quad = \int_0^T \!\!\int_{\Omega} v(x,t)\,\mathcal{L}u(x,t)\,dx\,dt 
		- \int_0^T \!\!\int_{\Omega} u(x,t)\,\mathcal{L}^* v(x,t)\,dx\,dt \\
		&\qquad\quad + \int_0^T \!\!\int_{\partial\Omega} \mathcal{T}^{cl+}_x u(x,t)\,v(x,t)\,d\sigma_x\,dt
		- \int_0^T \!\!\int_{\partial\Omega} \mathcal{T}^{cl+}_x v(x,t)\,u(x,t)\,d\sigma_x\,dt,
	\end{aligned}
\end{equation}
where the conormal derivative notation \eqref{Tcl}  was used in the boundary integrals over $\partial\Omega$.

We can now use the second Green identity \eqref{E4.13}, to deduce the parametrix-based third Green identity, that is a parametrix-based counterpart of the integral representation formula, similar to the variable-coefficient elliptic equations (see, e.g., \cite{mikhailov2002localized}, \cite{chkadua2009analysis2}).  
To this end,
for a fixed point $(y,\tau)\in Q$, we take
\begin{equation*}
	v(x,t) = P^*(x,y,t,\tau), 
\end{equation*}
Due to the singularity of the parametrix at the point $(x,t) = (y,\tau)$, we follow the argument of \cite[p.~19]{DohrSteinbachNiino+2019+1+60} and employ the second Green identity \eqref{E4.13} with $T=\tau - \varepsilon$, that is, in a truncated space--time cylinder
\begin{equation*}
	Q_{\tau-\varepsilon} := \Omega \times (0, \tau - \varepsilon),
\end{equation*}
where $0 < \varepsilon < \tau $.
This yields
	\begin{align}\label{limiting procedure of representation formula}
		&\int_{\Omega} u(x, \tau-\varepsilon)\, P^*(x,y,\tau-\varepsilon,\tau)\,dx 
		- \int_{\Omega} u(x,0)\, P^*(x,y,0,\tau)\,dx \nonumber \\
		&\quad = \int_0^{\tau-\varepsilon} \!\!\int_{\Omega} 
		\Big[ \mathcal{L}u(x,t)\, P^*(x,y,t,\tau)
		- R^*(x,y,t,\tau)\, u(x,t) \Big] dx\,dt  \nonumber\\
		&\qquad + \int_0^{\tau-\varepsilon} \!\!\int_{\partial\Omega} 
		\Big[ \mathcal{T}^{cl+}_x u(x,t)\, P^*(x,y,t,\tau)
		- \mathcal{T}^{cl+}_x P^*(x,y,t,\tau)\, u(x,t) \Big]
		d\sigma_x\,dt.
	\end{align}				
Taking the limit as $\varepsilon \to 0$, and using the same arguments as for the constant-coefficient heat equation 
(see, e.g., \cite[pp.~19--22]{DohrSteinbachNiino+2019+1+60}), the first term on the left-hand side of 
\eqref{limiting procedure of representation formula} converges to $u(y,\tau)$. 
Hence, passing to the limit in the remaining terms gives the following parametrix-based third Green identity for $(y,\tau) \in Q= \Omega \times (0,\infty)$:
	\begin{align}\label{representationformula}
		u(y,\tau)
		&= \int_{\Omega} P^*(x,y,0,\tau)\, u(x,0)\,dx
		+ \int_0^{\tau} \!\!\int_{\Omega} P^*(x,y,t,\tau)\,  \mathcal{L}u(x,t)\,dx\,dt \nonumber\\
		& + \int_0^{\tau} \!\!\int_{\partial\Omega}  P^*(x,y,t,\tau) \mathcal{T}^{cl+}_x u(x,t)d\sigma_x\,dt
		- \int_0^{\tau} \!\!\int_{\partial\Omega}\big(\mathcal{T}^{cl+}_x P^*\big)(x,y,t,\tau) u(x,t)	d\sigma_x\,dt\qquad \nonumber\\
		& - \int_0^{\tau} \!\!\int_{\Omega} R^*(x,y,t,\tau)\, u(x,t)\,dx\,dt.
	\end{align}	
Note that for $(y,\tau) \in Q$ and any $u \in C^{2}(\overline{Q})$, all integrals in \eqref{representationformula} are well defined.

\subsection{Parametrix-Based Potentials}
Let us introduce the space-time parametrix-based volume and layer potentials for linear non-stationary diffusion equations with variable coefficients and associated boundary operators.
For the constant-coefficient heat equation, the such potentials and boundary operators based on the fundamental solution are {\vl well known, see, e.g. \citep{costabel1990boundary, DohrSteinbachNiino+2019+1+60}}. 
{\bl Note that {\vl similar to \citep{costabel1990boundary, DohrSteinbachNiino+2019+1+60}, in definitions of the} parametrix-based potentials we employ a parametrix for the adjoint of the original differential operator under consideration.} 

The third Green identity \eqref{representationformula} can be re-written as
\begin{align}\label{3rdGreen}
		u(y,\tau)
		= \big(Mu(\cdot,0)\big)(y,\tau)
		&+ \big(\mathscr{P}  \mathcal{L}u\big)(y,\tau)
\nonumber\\
		&+ \big(V \mathcal{T}^{cl+}_x u\big)(y,\tau)
		- \big(W u\big)(y,\tau)
		- \big(\mathscr{R} u\big)(y,\tau).
	\end{align}	
Here {\bf $M$ is the parametrix-based initial potential operator} defined as
\begin{equation}\label{initial potential}\vl
			\left(M w_{0}\right)(y, \tau):=\int_{\Omega} P^*(x,y,0, \tau) w_{0}(x) \mathrm{d} x \quad \text { for }(y, \tau) \in \mathbb{R}^{d+1}.
		\end{equation}
Further, {\bf $\mathscr{P}$ is  the Newtonian and $\mathscr{R}$ is the remainder potential operators} defined, respectively, as
		\begin{equation}\label{Newtonian-lp}
			\left(\mathscr{P} w\right)(y,\tau):=\int_{0}^\tau\int_{\Omega} P^*(x,y,t,\tau) w(x, t) \,dx\,dt,
		\end{equation}
		\begin{equation}\label{Remainder-lp}
			\left(\mathscr{R} w\right)(y,\tau):=\int_{0}^\tau\int_{\Omega} R^*(x,y,t,\tau) w(x, t) \,dx\,dt
		\end{equation}
for $(y, \tau) \in \mathbb{R}^{d+1}$.
Finally, {\bf $V$ is the single layer and $W$ is the double layer operators}  defined, respectively,  
for $(y,\tau) \in \left(\mathbb{R}^{\MG d}\backslash \partial\Omega \right)\times\R$ as
		\begin{equation}\label{Single-lp}
			(V \widetilde w)(y, \tau) :=\int_{0}^\tau\int_{\partial\Omega} P^*(x,y,t,\tau) \widetilde w(x, t) \,d\sigma_x\, dt,
		\end{equation}
		\begin{equation}\label{Double-lp}
			(W \widetilde w)(y, \tau) :=\int_{0}^\tau\int_{\partial\Omega}\big(\mathcal{T}^{cl+}_xP^*\big)(x,y,t,\tau) \widetilde w(x, t) \,d\sigma_x\, dt.
		\end{equation}

		\subsection{Relations with Potentials of the Constant-Coefficient Heat Equation}
\label{sec 2.2}

Relationships between the {\vl operators} in \eqref{initial potential}–\eqref{Double-lp} and their constant-coefficient counterparts play a key role for exploring the properties of these operators. 
{\vl To this end, let us introduce the new variables}	
		\begin{equation}\label{change of variables}
			\tau_a:=\tau_{a(y,\tau)}:=a(y,\tau)\tau   \text{ and } t_a:=t_{a(y,\tau)}:=a(y,\tau)t.
		\end{equation}

		This implies that $dt=\frac{1}{a(y,\tau)} d t_a$, $a(y,\tau)(\tau - t)= \tau_a-t_a$, and 
$P^*(x,y,t,\tau)$ can be expressed in terms of fundamental solution of the  heat equation with the constant coefficient, $\underline{a}=1$, as
		\begin{multline}\label{P=P-h}
			P^*(x,y,t,\tau) =P_1^*(x,y, t_a, \tau_a)= 
\frac{\Theta(\tau_a-t_a)}{(4 \pi (\tau_a-t_a))^{d / 2}} \exp \left(\frac{-|y-x|^2}{4 (\tau_a-t_a)}\right)
\\
=P_1^*(x,y, t_a, a(y,\tau)\tau)= %\begin{cases}
\frac{\Theta(a(y,\tau)\tau-t_a)}{(4 \pi (a(y,\tau)\tau-t_a))^{d / 2}} \exp \left(\frac{-|y-x|^2}{4 (a(y,\tau)\tau-t_a)}\right).
		\end{multline}\label{E3.18}
		The {\vl integral potential densities} can be written as 
		\begin{align}\label{change of var-3}
			w(x, t)=\vl w(x,t_a/a(y,\tau)) =:w_a(x, t_a).
		\end{align}

Boundary of the corresponding parabolic cylinder $Q_a$ is denoted by $\Sigma_a=\partial\Omega\times (0,T_a)$. 
With these changes of variables, we can write {the parametrix-based} potentials \eqref{initial potential}--\eqref{Double-lp} in terms of the potentials $M_{\underline a}, V_{\underline a}, W_{\underline a}$ for the heat equation with a constant coefficient $\underline a$ (mostly with ${\underline a}=1$, and properties of such potentials can be found, e.g., in \citep{costabel1990boundary, DohrSteinbachNiino+2019+1+60}).	

		\begin{lemma}\label{potential relations Lemma}
Let functions $w_0$, $w$ and $\tilde w$ are defined, respectively on $\Omega$, $Q$ and $\Sigma$, and are sufficiently regular.
Then volume potentials satisfy the following  relations for $(y,\tau) \in \mathbb{R}^{d+1}$,
			\begin{align}
				{(M w_0)(y, \tau)}&={\left(M_1 w_0\right)\left(y, a(y,\tau)\tau\right)} 
; \label{initial pr}\\
				(\mathscr{P} w)(y, \tau)&=	\frac{1}{a(y,\tau)}\left(\mathscr{P}_1 w_a\right)\left(y, a(y,\tau)\tau\right)
; \label{NewtonianLPR}
			\end{align}			
\begin{align}\label{E3.22}
	(\mathscr{R}w)(y,\tau)=  
\sum_{i=1}^d\Big[\partial_{y_i}\mathscr{P}_{\underline a}\left(w\partial_i a\right)(y,\tau)\Big] _{\underline{a}=a(y,\tau)}
+(\widetilde{\mathscr R}w)(y,\tau),
\end{align}
where
\begin{align}\label{E3.23}
(\widetilde{\mathscr R}w)(y,\tau)&:=\int_{0}^\tau\int_{\Omega} \widetilde R^*(x,y,t,\tau) w(x, t) \,dx\,dt,
\end{align}
and $\widetilde R^*(x,y,t,\tau)$ is defined in \eqref{Reminder*~}.

Boundary layer potentials satisfy the following  relations for $(y, \tau)  \in\left(\mathbb{R}^{d}\backslash\partial\Omega\right) \times\R$,
			\begin{align}
				\left(V \widetilde w\right)(y,\tau)&=\frac{1}{a(y,\tau)} \left(V_1 \widetilde w_a\right)(y,a(y,\tau)\tau); 
\label{SingleLPR}\\
				\left(W\widetilde w\right)(y,\tau)&=\frac{1}{a(y,\tau)} W_1(a \widetilde w_a)(y,a(y,\tau)\tau).
\label{DoubleLPR}
			\end{align}	
		\end{lemma}	

		\begin{proof}
			 Relation  \eqref{initial pr} follows from 		change of variable in \eqref{P=P-h} when we  write $(Mw_0)(y,\tau)$ as follows:
				\begin{align*}
					(Mw_0)(y,\tau)&=\int_{\Omega}  P^*(x,y,0,\tau) w_0(x)\,d x\\
					&=\int_{\Omega}  P^*_1(x,y,0,a(y,\tau)\tau) w_0(x)\,d x\\
					&=(M_1w_0)(y,a(y,\tau)\tau).
				\end{align*}

For $\mathscr{P} w$, by \eqref{Newtonian-lp} and change of variables (see \eqref{change of variables} - \eqref{change of var-3}), we have
			\begin{equation*}
				\begin{split}
					\left(\mathscr{P} w\right)(y,\tau)&=\int_{0}^\tau\int_{\Omega} P^*(x,y,t,\tau) w(x, t) \,dx\,dt\\
					&=\frac{1}{a(y,\tau)}\int_{0}^{a(y,\tau)\tau}\int_{\Omega} P^*_1(x,y,t_a,a(y,\tau)\tau) w_a(x, t_a)  \,dx\,dt_a,\\
					&=\frac{1}{a(y,\tau)}\left(\mathscr{P}_1 w_a\right)(y,a(y,\tau)\tau).
				\end{split}
			\end{equation*}
			This proves relation \eqref{NewtonianLPR}. 

To prove relation \eqref{E3.22}, let us substitute $R^*(x,y,t,\tau)$ from \eqref{Reminder*} to \eqref{Remainder-lp} and obtain
\begin{align}\label{Reminder-op}
	(\mathscr{R}w)(y,\tau)&=\int_{0}^\tau\int_{\Omega} R^*(x,y,t,\tau) w(x, t) \,dx\,dt
\nonumber\\
&=\int_{0}^\tau\int_{\Omega}\left[\sum_{i=1}^{d} -\dfrac{\partial a(x,t)}{\partial x_i}\dfrac{\partial  P^*(x,y,t,\tau)}{\partial x_i}
+\widetilde R^*(x,y,t,\tau)\right]w(x, t) \,dx\,dt.
\end{align}	
Taking into account that by \eqref{E2.39} and \eqref{9.1*},
\begin{align*}
\dfrac{\partial  P^*(x,y,t,\tau)}{\partial x_i}
&=\dfrac{\partial}{\partial x_i}\bar P(x-y,\tau-t;a(y,\tau)) 
=\left[\dfrac{\partial}{\partial x_i}\bar P(x-y,\tau-t;\underline a)\right]_{\underline a=a(y,\tau)}
\\
&=-\left[\dfrac{\partial}{\partial y_i}\bar P(x-y,\tau-t;\underline a)\right]_{\underline a=a(y,\tau)}
=-\left[\dfrac{\partial}{\partial y_i} P_{\underline a}^*(x,y,t,\tau)\right]_{\underline a=a(y,\tau)}
\end{align*}
and employing notation \eqref{E3.23}, we obtain
\begin{align}\label{Reminder-op2}
	(\mathscr{R}w)(y,\tau)
=\left[\dfrac{\partial}{\partial y_i}\int_{0}^\tau\int_{\Omega}\sum_{i=1}^{d} P_{\underline a}^*(x,y,t,\tau)
w(x, t)\dfrac{\partial a(x,t)}{\partial x_i} \,dx\,dt\right]_{\underline a=a(y,\tau)}
+(\widetilde{\mathscr R}w)(y,\tau),
\end{align}	
leading to \eqref{E3.22}.
				
			Next, we prove relations \eqref{SingleLPR} and \eqref{DoubleLPR}.
			For any point $(y,\tau) \in (\mathbb{R}^{d} \setminus \partial\Omega) \R$, by combining  \eqref{Single-lp} and the change of variables formulas \eqref{change of variables}–\eqref{change of var-3}, we obtain
			\begin{equation*}
				\begin{split}
					(V\widetilde w)(y, \tau)
					&:= \int_{0}^\tau\int_{\partial\Omega} P^*(x,y,t,\tau) \widetilde w(x, t)\, d\sigma_x\, dt \\
					&= \int_{0}^{a(y,\tau)\tau}\int_{\partial\Omega} P_1^*(x,y, t_a, a(y,\tau)\tau) \widetilde w_a(x, t_a) \frac{1}{a(y,\tau)}\, d\sigma_x\, dt_a \\
					&= \frac{1}{a(y,\tau)} \int_{0}^{a(y,\tau)\tau}\int_{\partial\Omega} P_1^*(x,y, t_a, a(y,\tau)\tau) \widetilde w_a(x, t_a)\, d\sigma_x\, dt_a \\
					&= \frac{1}{a(y,\tau)} (V_1 \widetilde w_a)(y, a(y,\tau)\tau).
				\end{split}
			\end{equation*}

Similarly, taking into account \eqref{Tcl}, we obtain from \eqref{Double-lp} for $(y,\tau) \in (\mathbb{R}^{d} \setminus \partial\Omega)\times \R$,
			\begin{equation*}
				\begin{split}
					(W\widetilde w)(y, \tau)
					&:= \int_{0}^{\tau}\int_{\partial\Omega} \mathcal{T}^+_{x} P^*(x,y,t,\tau) \widetilde w(x, t)\, d\sigma_x\, dt \\
					&= \int_{0}^{\tau}\int_{\partial\Omega} \sum_{i=1}^{d} a(x,t) n_i(x)
					\frac{\partial P^*(x,y,t,\tau)}{\partial x_i} \widetilde w(x, t)\, d\sigma_x\, dt \\
					&= \int_{0}^{a(y,\tau)\tau}\int_{\partial\Omega} \sum_{i=1}^{d} n_i(x)
					\frac{\partial P_1^*(x,y, t_a, a(y,\tau)\tau)}{\partial x_i}  a(x,t) \widetilde w_a(x, t_a)
					\frac{1}{a(y,\tau)}\, d\sigma_x\, dt_a \\
					&= \frac{1}{a(y,\tau)} W_1(a \widetilde w_a)(y, a(y,\tau)\tau).
				\end{split}
			\end{equation*}
		\end{proof} 

			As in the study of boundary-domain integral equations for linear elliptic PDEs with variable coefficients, the layer potential relations in Lemma \ref{potential relations Lemma} allow one to analyse the mapping properties of integral potentials for linear non-stationary diffusion problems with variable coefficients via the corresponding integral potentials for the heat equation with constant coefficients.

		\section{Conclusion}
		This paper has established basic elements for a space-time boundary-domain integral equation framework for linear non-stationary heat equations with a variable coefficient. Our primary aim is the construction of an explicit space-time parametrix and the proof that the corresponding remainder operator is at most weakly singular. 
		
		Furthermore, we have demonstrated 
explicit relations between the volume and layer potentials for the variable-coefficient heat equation and their well-known counterparts for the constant-coefficient heat equation. These include the initial, Newtonian, single-layer, and double-layer
potentials.
		
		The results presented herein provide the essential theoretical basis for 
the {\MG next} parts of the research.
{\MG They} should include the detailed analysis of the mapping properties of these parametrix-based potentials in anisotropic Sobolev spaces, and the subsequent proof of the equivalence between a given initial-boundary value problem and the derived systems of space-time boundary-domain integral equations, leading to analysis of the systems unique solvability for Dirichlet, Neumann, and mixed boundary conditions.

\paragraph{Acknowledgment:} The second author gratefully acknowledges the support from The QJMAM Fund for Applied Mathematics allowing him to present some of the paper results at the IMSE 2026 conference in Matera, Italy, and {\MG to} finalise the paper.

	\end{document}